\documentclass[journal]{IEEEtran}

\usepackage{epsfig,amsfonts,color}
\usepackage[final]{pdfpages}
\usepackage{amsmath}
\usepackage{blkarray}
\usepackage{bm}
\usepackage{algorithm}
\usepackage{algpseudocode}
\usepackage{multirow}
\renewcommand{\algorithmicfunction}{\textbf{Function}}
\algdef{SE}[TASK]{Task}{EndTask}%
   [2]{\algorithmicfunction\ \tt{#1}\ifthenelse{\equal{#2}{}}{}{(#2)}}%
   {\algorithmicend\ \algorithmicfunction}%
\usepackage{graphicx}
\usepackage[belowskip=-10pt]{caption}
\usepackage{subcaption}
\usepackage{placeins}

\usepackage{amssymb,url}
\usepackage[colorlinks=true,linkcolor=blue,citecolor=blue,urlcolor=blue]{hyperref}

\usepackage{xcolor}
\usepackage{multirow}
\usepackage{mathtools}
\usepackage{fancyhdr}

\newcommand{\be}{\begin{equation}}
\newcommand{\ee}{\end{equation}}
\newcommand{\bea}{\begin{eqnarray}}
\newcommand{\eea}{\end{eqnarray}}

\newcommand{\pref}[1]{(\ref{#1})}
\newcommand{\bvec}{\left(\begin{array}{c}}
\newcommand{\evec}{\end{array}\right)}
\newcommand{\bsub}{\begin{subequations}}
\newcommand{\esub}{\end{subequations}}

\newcommand{\bx}{\boldsymbol{x}}
\newcommand{\by}{\boldsymbol{y}}
\newcommand{\bz}{\boldsymbol{z}}
\newcommand{\bb}{\boldsymbol{b}}
\newcommand{\bc}{\boldsymbol{c}}
\newcommand{\bpi}{\boldsymbol{\pi}}
\newcommand{\blambda}{\boldsymbol{\lambda}}
\newcommand{\bmu}{\boldsymbol{\mu}}

\newcommand{\br}{\boldsymbol{r}}
\newcommand{\ml}{\mathcal{L}}

\usepackage{accsupp}

\usepackage{listings}
\DeclareCaptionFont{white}{\color{white}}
\DeclareCaptionFormat{listing}{{\parbox{\dimexpr\textwidth-2\fboxsep\relax}{#1#2#3}}}

\usepackage[dvipsnames]{xcolor}

\lstdefinelanguage{julia}
{
alsoletter={!},
keywordsprefix=\@,
morekeywords={
exit,whos,edit,load,is,isa,isequal,typeof,tuple,ntuple,uid,hash,finalizer,convert,promote,
subtype,typemin,typemax,realmin,realmax,sizeof,eps,promote_type,method_exists,applicable,
invoke,dlopen,dlsym,system,error,throw,assert,new,Inf,Nan,pi,im,begin,while,for,in,return,
break,continue,macro,quote,let,if,elseif,else,try,catch,end,bitstype,ccall,do,using,module,
import,export,importall,baremodule,immutable,local,global,const,Bool,Int,Int8,Int16,Int32,
Int64,Uint,Uint8,Uint16,Uint32,Uint64,Float32,Float64,Complex64,Complex128,Any,Nothing,None,
function,type,typealias,abstract,get_node,add_edge,create_estimation_model,set_solution,
solve, get_solution, solve_ss_problem, create_estimation_problem, addnode, Partition, 
assemble_optigraph, local_subgraphs, apply_partition, OptiGraph, optimizer_with_attributes,
BendersAlgorithm, optimize, aggregate, source_graph, aggregate_to_depth, add_subgraph, fill, set_to_node_objectives, set_optimizer, RemoteOptiGraph, run_algorithm!, setdiff, all_variables, Symbol, fetch, load_operation_optinode_remote!, name, local_graph, GenericAffExpr, operation_model!, replace, intersect, findfirst, Dict
},
morekeywords = [2]{triggered_by,compute_time,trigger_during_busy,send_on,delay,send_wait,start},
sensitive=true,
morecomment=[l]{\#},
morestring=[b]',
morestring=[b]"
}

\lstdefinelanguage{julia}
{
alsoletter={!},
keywordsprefix=\@,
morekeywords={
exit,whos,edit,load,is,isa,isequal,typeof,tuple,ntuple,uid,hash,finalizer,convert,promote,
subtype,typemin,typemax,realmin,realmax,sizeof,eps,promote_type,method_exists,applicable,
invoke,dlopen,dlsym,system,error,throw,assert,new,Inf,Nan,pi,im,begin,while,for,in,return,
break,continue,macro,quote,let,if,elseif,else,try,catch,end,bitstype,ccall,do,using,module,
import,export,importall,baremodule,immutable,local,global,const,Bool,Int,Int8,Int16,Int32,
Int64,Uint,Uint8,Uint16,Uint32,Uint64,Float32,Float64,Complex64,Complex128,Any,Nothing,None,
function,type,typealias,abstract,get_node,add_edge,create_estimation_model,set_solution,
solve, get_solution, solve_ss_problem, create_estimation_problem, addnode, Partition, 
assemble_optigraph, local_subgraphs, apply_partition, OptiGraph, optimizer_with_attributes,
BendersAlgorithm, optimize, aggregate, source_graph, aggregate_to_depth, add_subgraph, fill, set_to_node_objectives, set_optimizer, RemoteOptiGraph, run_algorithm!, setdiff, all_variables, Symbol, fetch, load_operation_optinode_remote!, name, local_graph, GenericAffExpr, operation_model!, replace, intersect, findfirst, Dict
},
morekeywords = [2]{triggered_by,compute_time,trigger_during_busy,send_on,delay,send_wait,start},
sensitive=true,
morecomment=[l]{\#},
morestring=[b]',
morestring=[b]"
}

\usepackage{authblk}
\usepackage{amsthm}
\theoremstyle{plain}

\theoremstyle{definition}

\usepackage{orcidlink}

\usepackage{titlesec}
\titleformat{\subsubsection}
[runin]
{\normalfont\itshape}
{\arabic{subsubsection})}
{2pt}{}[:\hspace*{2pt}]
\titlespacing*{\section}{0pt}{*0.9}{*0.9}
\titlespacing{\subsection}{0pt}{*0.9}{*0.9}
\titlespacing{\subsubsection}{0pt}{*0.9}{*0.9}
\begin{document}

\title{\bf \Large Analyzing Performance and Scalability of Benders Decomposition for Generation and Transmission Expansion Planning Models}

\author{David L. Cole\orcidlink{0000-0002-5703-6958},
Michael Lau\orcidlink{0000-0002-2212-4069},
Xinliang Dai\orcidlink{0000-0001-8669-7942},
Sambuddha Chakrabarti\orcidlink{0000-0002-8916-5076},
Jesse D. Jenkins\orcidlink{0000-0002-9670-7793}\thanks{This work was supported by a grant from Clean Grid Initiative and Breakthrough Energy\\
David Cole, Sambuddha Chakrabarti, Xinliang Dai, and Jesse Jenkins are with the Andlinger Center for Energy and the Environment, Princeton University, Princeton, NJ 08540 USA.\\
Michael Lau and Jesse Jenkins are with the Mechanical and Aerospace Engineering Department, Princeton University, Princton, NJ 08540 USA.},
}


\maketitle

\begin{abstract}
Generation and Transmission Expansion Planning (GTEP) problems co-optimize generation and transmission expansion, enabling them to provide better planning decisions than traditional Generation Expansion Planning  (GEP) or Transmission Expansion Planning (TEP) problems. However, GTEPs are computationally complex because they must represent network structure and power flow constraints with higher resolution than GEPs while considering generation expansion, unlike TEPs. As a result, GTEPs can be intractable without decomposition. Benders Decomposition (BD) has been applied to expansion planning problems, with various methods applied to accelerate convergence. In this work, we test strategies for improving the performance of BD on GTEP models with nodal resolution and DC optimal power flow (DCOPF) constraints, and present an alternative approach for handling the bilinear constraints that can result in these problems. These tests included combinations of using generalized Benders decomposition (GBD), hot-starting via a transport constrained model, using linear relaxations of the master problem, and using regularization. We test these methods on mixed integer linear programming GTEP models with 24, 73, and 146 bus systems with sixteen representative weeks (up to 10 million continuous variables and 400 mixed-integer decisions) which we find to be computationally intractable without decomposition. With selected accelerated Benders decomposition approaches, the problems can be solved to under a 1\% gap in as little as 5 hours. We find that GBD can improve time to convergence, and that for both BD and GBD, combining hot-starting from a transport model and linear relaxations of the master problem significantly improved performance and were necessary for recovering high quality solutions. The results also suggest that using regularization on these initial hot-starting and relaxation steps and then turning it off after they are complete was generally the best combination of strategies for both BD and GBD. 
\end{abstract}


\section{Introduction}

Generation and transmission expansion planning (GETP) models can identify plans with superior financial and environmental improvement relative to traditional modeling approaches that consider generation expansion planning (GEP) or transmission expansion planning (TEP) models alone. Traditional methods generally take a reactive approach, where transmission expansion is done after generation expansion is solved, or an iterative approach, where GEP and TEP models are iteratively solved. In contrast, GTEP models seek to simultaneously co-optimize new generation and transmission asset investment and retirement decisions in a future grid over a given time horizon \cite{gonzalez2020review, krishnan2016co}.
As new generation siting and sizing decisions can impact the location of new transmission capacity and vice versa, co-optimizing generation and transmission expansion can identify optimal expansion plans that are billions of dollars less expensive than isolated planning approaches \cite{krishnan2016co, mehrtash2023necessity, spyrou2017benefits}. The results of such GTEP models can then be used to inform grid transition strategies, assess new or emerging technologies, evaluate potential policies, or value asset siting and sizing decisions.


Despite their value, GTEPs are typically computationally infeasible to solve because they must represent systems with sufficient network resolution and  constraint complexity to reflect realistic physical power flow behavior yet maintain sufficient temporal detail to accurately evaluate generation investment decisions. Significant work has been done to produce informative models that simplify one or more of these domains (e.g., using  representative days, pre-filtering new line options, relaxing the underlying formulation), but transmission expansion decisions rely on accurately capturing power flow dynamics. While many GEP problems use a transport model (``pipe-flow''; i.e., power flows on lines are only constrained by a power balance and upper and lower capacity limits) for transmission constraints, GTEP models require more accurate power flow dynamics by modeling AC power flow or a DC optimal power flow linearization. When considering discrete transmission build decisions, capturing these operations can also result in bilinear constraints that complicate solution. Further, these power flow dynamics occur at the nodal (i.e., bus) level, so high spatial resolutions can be important for GTEP. Consequently, improving mathematical solution abilities for GTEP problems with optimal power flow and high spatial resolution is valuable to better understand and construct reliable and affordable future power systems.

Decomposition algorithms have been a powerful tool to push the bounds of existing solution capabilities in GTEP problems. Progressive hedging (PH) \cite{rockafellar1991scenarios} is used in stochastic studies, where first-stage variables are determined and then fixed in scenario-specific subproblems. PH is an iterative algorithm and has been used successfully in GTEP studies \cite{mahdavi2018transmission,zuluaga2024parallel,zuluaga2025nodal}. Benders Decomposition (BD) \cite{benders1962partitioning} is an iterative cutting plane algorithm that separates the optimization problem into a master problem and one or more subproblems. The master and subproblems are iteratively solved, passing primal and dual information between them and adding cuts to the master problem until a termination criterion is reached. BD has been used in many GTEP studies over the decades \cite{cho2022recent, emdadi2025benders, jenabi2015acceleration, li2022mixed, najjar2018coordinated, pereira1985decomposition, roh2009market, sasanpour2025accelerating, thome2013decomposition}
with varying formulations in network representations, system sizes, and technologies considered. Other studies have used BD for storage siting and transmission expansion \cite{de2022block,macrae2016benders, wu2024high, wu2025high}. Algorithmic advances have been made to accelerate solution times, such as using simplified models to generate initial cuts \cite{jenabi2015acceleration,romero2002hierarchical, thome2013decomposition} or linear programming (LP) relaxations \cite{lumbreras2016solve}. While Lumbreras et al. \cite{lumbreras2016solve} considered the impacts of various speedup strategies on transmission-only expansion problems, to our knowledge, no study has compared the performance of various speedup strategies for BD on GTEP models, which feature a larger number of expansion decision variables.  

In this work, we evaluate four different speedup strategies for BD with DCOPF-constrained GTEP models with discrete decision variables for new lines. The systems we solved consider discrete generator siting and sizing decisions, as well as the ability to reconductor a subset of the existing lines to increase capacity. We model operational decisions using 16 representative weeks to reduce the temporal domain, and the methods are applied to three different system sizes of 24, 73, and 146 buses based on the Reliability Test System (RTS) \cite{barrows2019ieee}. The first strategy is to use GBD to solve the subproblems with quadratic constraints. Because BD is grounded in linear duality theory, it is not applicable directly to the nonlinear case. We propose using GBD to solve these problems as it can potentially produce stronger cuts by maintaining the nonlinear constraint in the subproblems. We also test three additional speedup strategies: i) using a transport model to ``hot-start'' the DCOPF model, ii) using an LP relaxation of the master problem to generate initial cuts for the DCOPF model, and iii) using a regularization scheme to generate better master problem decisions. Although one or more of these speedup approaches have been used in individual studies, to our knowledge, they have not been directly compared for GTEP models. Furthermore, to our knowledge, none have been applied to GBD in the literature.

This work differentiates itself from the literature in several ways. While Kim et al. \cite{kim2015integrated} suggest the idea of using GBD for DCOPF-constrained GTEP problems, they do not provide information on their implementation or numerical results. We report computational results for a much larger system than they consider, and do so in a replicable way for other researchers. Several works in TEP have used approaches that were effectively GBD, though they do not use that name. Alizadeh-Mousavi et al. \cite{alizadeh2016efficient}, Orfanos et al. \cite{orfanos2012transmission}, and Zhan et al. \cite{zhan2017fast} each apply BD and use a different formulation to recover cuts. Alizadeh-Mousavi et al. \cite{alizadeh2016efficient} and Orfanos et al. \cite{orfanos2012transmission} use sensitivity information to generate cuts, but do not give proofs that these cuts are identical to the classic cut obtained by GBD for their system. Zhan et al. \cite{zhan2017fast} use an LP relaxation of the nonlinear subproblem to recover a Lagrangian multiplier. In Gu et al. \cite{gu2012coordinating}, they likewise use a sensitivity metric, but they apply it to a GTEP problem that considered ancillary services and iterated between a generation expansion problem and transmission expansion problem to reach results rather than co-optimizing in the same optimization problem, and they do so on a smaller system than what we will show in this work. While considering the nonlinear subproblem, none of these works prove that their resulting cuts are a lower bound on the problem. Because the nonlinear form of the subproblem is nonconvex (it involves a  quadratic equality constraint), GBD may not converge to a global or even local optimum \cite{bagajewicz1991generalized, sahinidis1991convergence}. In our work, we do provide a lower bound derived from a convex relaxation featuring transport constraints for power flows. Lastly, we note the work of Romero and Monticelli \cite{romero2002hierarchical} who introduced a hierarchical approach for a model with TEP. They note the nonlinear form of the original problem and solve a transport model to generate initial cuts and then solve the DCOPF problem using a linear form of the problem and sensitivity information to generate cuts. We expand on these previous works by applying GBD directly to a GTEP problem. We provide a guaranteed lower bound by using the solution of a transport model within the algorithm. Finally, while other works have used speedup strategies for BD on GTEP models, like using the transport model for hot-starting or using an LP relaxation or regularization, none have provided a comparative analysis of the different strategies as we do.

Other researchers have successfully solved larger systems (in terms of number of buses) and produced promising results in the area of GTEP and TEP, and we also note the placement of our work among these others. Wu et al. \cite{wu2024high, wu2025high} solved a 2,000 bus transmission expansion case with storage siting using BD. In their studies, they consider storage and transmission expansion (not new generation), and their expansion of new lines is modeled as discrete increased capacity on existing lines rather than creating new lines. In our work, we consider new lines and new generating capacity in addition to reconductoring of existing lines. Zuluaga et al. \cite{zuluaga2024parallel} and Musselman et al. \cite{musselman2025climate} solve over 8,000 buses, the largest cases of which we are aware in terms of the number of buses. They include new line construction but use a transport model for power flow rather than DCOPF. In a more recent study by Zuluaga et al. \cite{zuluaga2025nodal}, they do consider DCOPF constraints on a smaller system with new line construction. In these three latter studies, they consider a stochastic model and use representative days. Further, they use progressive hedging to solve their models rather than BD. Their results suggest that the systems they solve are not scalable without decomposition schemes like PH or BD. In addition, these studies use large compute resources (e.g., more than 3000 CPU cores for the largest instance in \cite{zuluaga2024parallel}) to solve their models with varying success in terms of optimality gap (ranging from 2\% to over 10\%). In comparison, the system we consider has fewer buses but considers representative weeks of operational decisions rather than days, and we solve our system to a smaller optimality gap.

In short, our paper advances the literature by i) analyzing the empirical performance of different speedup strategies for BD on GTEP problems to ensure scalable solution and convergence, and  ii) analyzing the performance of GBD on GTEP problems. Our results suggest that hybrid combinations of hot-starting with the transport model, using an LP relaxation of the master problem, and regularization can greatly improve performance of BD and GBD and that GBD can provide improved performance over BD in some cases. 

This paper is structured as follows. Section \ref{sec:problem_formulation} gives the problem formulation. Section \ref{sec:GBD} introduces GBD and its implementation in our code. Section \ref{sec:case_study} presents case studies with three different system sizes to compare performance of GBD with BD and discusses the results. Section \ref{sec:conclusion} gives conclusions and future work.

\section{Problem Formulation}\label{sec:problem_formulation}

We consider a GTEP problem with investment and operational decisions over a period of time. The time period is split up into a set of sub-periods, $W$, featuring operational decisions at hourly resolution. We denote the vectors of new generation investment decisions as $\by$, of new transmission line investment decisions as $\bz$, and of reconductoring decisions as $\br$. The vector of operational decisions for a subperiod $w \in W$ is given as $\bx_w$, which includes dispatch decisions, power flows across transmission lines, bus phase angles, unit commitment decisions, and slack variables (e.g., for load shedding). The compact mathematical formulation for this problem is given by \eqref{eq:CEM}.
\begin{subequations}\label{eq:CEM}
\begin{flalign}
    \qquad  \min &\; \bc^\top_y \by + \bc^\top_z \bz + \bc^\top_r \br + \sum_{w \in W} \beta_{w}\bc^\top_w \bx_w \label{eq:CEM_obj}\\
    \textrm{s.t.} &\; D_w \bx_w \le \boldsymbol{d}_w, & w \in W \label{eq:CEM_operation_cons} \\
    &\; A_w \bx_w + B_w \by \le \bb_w, & w \in W \label{eq:CEM_generator_linking}\\
    &\; g_w(\bx_w, \bz, \br) = 0, & w \in W \label{eq:CEM_NLP_equalities} \\
    &\; h_w(\bx_w, \bz, \br) \le 0, & w \in W \label{eq:CEM_NLP_inequalities} \\
    &\; \underline{\bx} \le \bx_w \le \overline{\bx}, & w \in W \label{eq:CEM_op_vars_limits}\\
    &\; 0 \le \br \le \overline{\br} \label{eq:CEM_r_var_limits} \\
    &\; \by \in \mathbb{Z}_+^{n_y},   \label{eq:CEM_int_vars}\\
    &\; \bz \in \{0,1\}^{n_z}. \label{eq:CEM_bin_vars}
\end{flalign}
\end{subequations}
\noindent Here, \eqref{eq:CEM_obj} is the objective, which seeks to minimize the investment and operational costs, where $\bc_y$, $\bc_z$, $\bc_r$, and $\bc_w$ are costs and $\beta_w$ is a weighting on the operational sub-period $w$ to ensure operational costs are equivalent to a full-year of operations. \eqref{eq:CEM_operation_cons} are operational constraints for sub-period $w$, such as unit commitment and ramping for thermal generators. \eqref{eq:CEM_generator_linking} are linking constraints between new generator capacity and operational variables, such as limits on power production based on the size of the new generators. \eqref{eq:CEM_NLP_equalities} and \eqref{eq:CEM_NLP_inequalities} are constraints on the transmission operation that will be discussed more below, but which can be nonlinear in the DCOPF constrained case. \eqref{eq:CEM_op_vars_limits} and \eqref{eq:CEM_r_var_limits} are bounds on their respective variables (with $\underline{\cdot}$ and $\overline{\cdot}$ being vectors of lower or upper bounds). \eqref{eq:CEM_int_vars} and \eqref{eq:CEM_bin_vars} require their respective variables to be integer or binary, with $n_y$ and $n_x$ being the length of their respective vectors. The sub-periods in this formulation can be representative periods of varying lengths. Importantly, there are no linking constraints between sub-periods, only between the investment decisions and sub-periods. See \cite{pecci2025regularized} for more on how constraints linking sub-periods (e.g., emissions budgets or long-duration storage levels) can be resolved into the generalized form above by implementing budgeting decisions in the master problem.

The transmission constraints \eqref{eq:CEM_NLP_equalities} and \eqref{eq:CEM_NLP_inequalities} can have several forms. In this work, we will use DCOPF constraints and support multiple types of transmission expansion decisions. We denote the set of all corridors in the system as $\mathcal{C}$, where a corridor is a set of two buses, $(i,j)$, between which a transmission line is constructed or could be constructed. Corridors can support parallel lines, such that we will index variables by a line index rather than by a corridor. We will refer to potential new lines that can be constructed as ``candidates''. The set of all lines (existing and candidates) is given by $\ml$, with $\ml(i,j)$ representing the subset that are on corridor $(i,j)$ (we will maintain this notation for the other subsets of lines). The set of all candidate lines is $\ml_{cand}$, while the set of all existing lines which can be reconductored are $\ml_{rec}$. Reconductoring is a process of upgrading power lines (such as replacing with new cables) which increases the capacity of existing lines, and this is modeled by a continuous variable that can increase capacity by some limited amount (see also \eqref{eq:CEM_r_var_limits}). To maintain greater generality, we consider two different ways of adding candidate lines: i) A new line on a new or existing corridor which could be parallel to existing lines, in which case existing lines are represented as $\ml_{fixed}$ , and ii) candidate lines that would replace an existing line if it was constructed, in which case existing lines are represented as $\ml_{rep}$. 

The more detailed transmission constraints for \eqref{eq:CEM_NLP_equalities} and \eqref{eq:CEM_NLP_inequalities} for a time step $t$ of a sub-period $w$ (sub-period index will be omitted for simplicity) is given by \eqref{eq:DCOPF_bilinear}.
\begin{subequations}\label{eq:DCOPF_bilinear}
\begin{flalign}
    & D_{i, t} = \sum_{g \in \mathcal{R}_i} p_{g, t} - \sum_{(j,k) \in \mathcal{C}| j = i} \Bigg( \sum_{\ell \in \mathcal{L}(j,k)} f_{\ell,t}\Bigg) \notag\\&\qquad\qquad+ \sum_{(k,j) \in \mathcal{C}| j = i} \Bigg( \sum_{\ell \in \mathcal{L}(k,j)} f_{\ell, t}\Bigg) + s_t^+, \; i \in \mathcal{N} \label{eq:DCOPF_nodal_balance} \\
    & f_{\ell,t} = B_\ell \Delta\theta_{ij, t}, \qquad \qquad \;\,(i, j) \in \mathcal{C}, \ell \in \mathcal{L}_{fixed}(i,j) \label{eq:DCOPF_exist} \\
    & f_{\ell, t} = B_\ell \Delta \theta_{ij, t}z_{\ell}, \qquad \qquad (i, j) \in \mathcal{C}, \ell \in \mathcal{L}_{cand}(i,j) \label{eq:DCOPF_bilinear_cand} \\
    & f_{\ell,t} = B_\ell \Delta\theta_{ij, t}(1-z_{\psi(\ell)}), \quad (i,j) \in \mathcal{C}, \ell \in \mathcal{L}_{rep}(i,j) \label{eq:DCOPF_bilinear_replace} \\
    & \underline{f}_\ell \hspace{-1pt}- \hspace{-1pt}r_{\ell}^{low} \hspace{-1pt}-\hspace{-1pt} r_{\ell}^{high} \le f_{\ell,t} \le \overline{f}_\ell \hspace{-1pt}+ \hspace{-1pt}r_{\ell}^{low}\hspace{-1pt} +\hspace{-1pt} r^{high}_{\ell}, \, \ell \in \mathcal{L}_{rec} \label{eq:DCOPF_limits_reconductor} \\
    & \underline{f}_\ell \le f_{\ell,t} \le \overline{f}_\ell, \qquad \qquad \qquad \qquad \qquad \;\, \ell \in \mathcal{L} \setminus\mathcal{L}_{rec} \label{eq:DCOPF_limits_noreconductor}\\
    & -\frac{\pi}{4} \le \theta_{i,t} \le \frac{\pi}{4} \label{eq:DCOPF_angle_limits}
\end{flalign}
\end{subequations}
\noindent Here, \eqref{eq:DCOPF_nodal_balance} is the nodal power balance, where $\mathcal{N}$ is the set of buses in the system, $D_{i,t}$ is the power demand, $\mathcal{R}_i$ is the set of generation resources on bus $i$, $p_{g,t}$ is the power injection from resource $g$, $f_{\ell, t}$ is the power flow across line $\ell$, and $s^+_t$ are slack variables for load shedding. Slacks $s^+_t$are heavily penalized in the objective \eqref{eq:CEM_obj}. \eqref{eq:DCOPF_exist} are DCOPF constraints dictating power flow for existing (fixed) lines, where $B_\ell$ is the susceptance of line $\ell$ and $\Delta \theta_{ij,t}$ is the bus phase angle difference between buses $i$ and $j$. \eqref{eq:DCOPF_bilinear_cand} are power flow constraints for candidate lines, where $z_\ell$ is a binary build decision for line $\ell$, such that $f_{\ell,t}$ is zero if the line is not built ($z_\ell = 0$) and $B_\ell \Delta \theta_{ij,t}$ otherwise. \eqref{eq:DCOPF_bilinear_replace} are power flow constraints for lines replaceable by a candidate line, where $\psi(\ell)$ is a mapping to the candidate line which replaces $\ell$ if constructed (i.e., if $z_{\psi(\ell)} = 1$, then $f_{\ell, t} = 0$). Note that \eqref{eq:DCOPF_bilinear_cand} and \eqref{eq:DCOPF_bilinear_replace} are both bilinear terms, wherein the capacity decision $z$ is multiplied by the voltage angle decision $\theta_{ij,t}$. \eqref{eq:DCOPF_limits_reconductor} are upper and lower limits on power flow for reconductored lines, where there are two different reconductoring variables. These variables have different costs associated with them, effectively creating a piecewise cost function on the reconductoring. \eqref{eq:DCOPF_limits_noreconductor} are upper and lower bounds on lines that cannot be reconductored.  \eqref{eq:DCOPF_angle_limits} are limits on phase angles.

In many studies, the bilinear terms of \eqref{eq:DCOPF_bilinear_cand} and \eqref{eq:DCOPF_bilinear_replace} are linearized through use of a ``big M'' relaxation. The form of this relaxation for these constraints is given by \eqref{eq:DCOPF_linear}
\begin{subequations}\label{eq:DCOPF_linear}
\begin{align}
    &\; -M_\ell(1 - z_{\ell}) \le f_{\ell,t} - B_\ell \Delta \theta_{ij, t} \le M_\ell(1 - z_{\ell}), \notag\\
    &\qquad\qquad \qquad\qquad\qquad\qquad \, (i,j) \in \mathcal{C},  \ell \in \mathcal{L}_{cand}(i,j) \label{eq:DCOPF_bigM_cand1} \\
    &\; \underline{f}_\ell z_{\ell} \le f_{\ell, t} \le \overline{f}z_{t,\ell}, \quad \qquad\;\,(i,j) \in \mathcal{C},  \ell \in \mathcal{L}_{cand}(i,j) \label{eq:DCOPF_bigM_cand2} \\
    &\; -M_\ell z_{\psi(\ell)} \le f_{\ell,t} - B_\ell \Delta \theta_{ij, t} \le M_\ell z_{\psi(\ell)},\notag\\
    &\qquad\qquad \qquad\qquad\qquad\qquad \, (i,j) \in \mathcal{C},  \ell \in \mathcal{L}_{rep}(i,j) \label{eq:DCOPF_bigM_rep1} \\
    &\; \underline{f}_\ell (1 - z_{\psi(\ell)}) \le f_{\ell,t} \le \overline{f}_{\ell}(1-z_{\psi(\ell)}), \notag\\
    &\qquad\qquad \qquad\qquad\qquad\qquad \, (i,j) \in \mathcal{C},  \ell \in \mathcal{L}_{rep}(i,j) \label{eq:DCOPF_bigM_rep2}
\end{align}
\end{subequations}
\noindent Here, \eqref{eq:DCOPF_bigM_cand1} and \eqref{eq:DCOPF_bigM_cand2} replace constraint \eqref{eq:DCOPF_bilinear_cand}. When $z_\ell = 0$, $f_{\ell,t}$ is not constrained by the value of $\Delta \theta_{ij,t}$ in \eqref{eq:DCOPF_bigM_cand1} and is forced to equal zero by \eqref{eq:DCOPF_bigM_cand2}. When $z_{\ell} = 1$, $f_{\ell, t}$ is forced to equal $B_\ell \Delta \theta_{ij,t}$, with line limits set by \eqref{eq:DCOPF_bigM_cand2}. \eqref{eq:DCOPF_bigM_rep1} and \eqref{eq:DCOPF_bigM_rep2} replace \eqref{eq:DCOPF_bilinear_replace} by similar logic. Because of the bilinear term, solving \eqref{eq:DCOPF_bilinear} can require mixed-integer nonlinear solvers, whereas \eqref{eq:DCOPF_linear} can be solved with normal mixed-integer LP solvers, making the latter more efficient when solving \eqref{eq:CEM} without decomposition. Note that the value of $M_\ell$ can be important for numerical reasons \cite{binato2002new}. 

\section{Benders Decomposition} \label{sec:GBD}

In this section, we introduce conventional BD and then define the more generic GBD for nonlinear problems. We also introduce strategies for improving convergence of BD and GBD

\subsection{Base Formulation and Algorithm}

BD \cite{benders1962partitioning} is an iterative algorithm that decomposes a problem into a master problem and one or more subproblems. The master problem contains a set of complicating variables that, once fixed in the subproblems, make the subproblems easier to solve and separable. At each iteration, the master problem is solved, primal solutions of the master problem enter the subproblem(s) as fixed parameters, and the subproblem(s) are solved. Then, information (typically dual information) is passed back to the master problem to form cutting planes, which provide a lower bound approximation of the subproblem objective as a function of master problem decisions. This process is repeated until upper and lower bound estimates of the full problem objective converge to an acceptable gap. For our problems, we place the investment decisions ($\by$, $\bz$, and $\br$) in the master problem and fix their solutions in the subproblem.

Given a set of investment decisions from the master problem---$\by^{(k)}$, $\bz^{(k)}$, and $\br^{(k)}$---at iteration $k$, each operational sub-period takes the form given by \eqref{eq:CEM_subproblem}.
\begin{subequations}\label{eq:CEM_subproblem}
\begin{flalign}
    \qquad\quad \; \phi^{(k)}_w = \min &\; \beta_w \bc^\top_w \bx_w \\
    \textrm{s.t.} &\; D_w \bx_w \le \boldsymbol{d}_w, & w \in W  \\
    &\; A_w \bx_w + B_w \by \le \bb_w, & w \in W \\
    &\; g_w(\bx_w, \bz, \br) = 0, & w \in W \label{eq:CEM_subproblem_NLP_equality} \\
    &\; h_w(\bx_w, \bz, \br) \le 0, & w \in W \label{eq:CEM_subproblem_NLP_inequality} \\
    &\; \by = \by^{(k)} &(\bpi_w^{(k)}) \\
    &\; \bz = \bz^{(k)} &(\blambda_w^{(k)}) \\
    &\; \br = \br^{(k)} &(\bmu_w^{(k)}) \\
    &\; 0 \le \bx_w \le \overline{\bx}, & w \in W
\end{flalign}
\end{subequations}
\noindent Here, $\phi^{(k)}_w$ is the objective value of subproblem $w$ at iteration $k$ and $\bpi^{(k)}_w$, $\blambda^{(k)}_w$, and $\bmu^{(k)}_w$ are dual variables associated with their respective constraints. Constraints \eqref{eq:CEM_subproblem_NLP_equality} and \eqref{eq:CEM_subproblem_NLP_inequality} take the form of \eqref{eq:DCOPF_linear}. Because there is no linking between the operational sub-periods $w \in W$, each subproblem can be solved independently of and in parallel with the other subproblems \cite{Jacobson2024AConstraints}. 

The master problem takes the form given by \eqref{eq:CEM_master}.
\begin{subequations}\label{eq:CEM_master}
\begin{align}
    L^{(k)} = \min &\; \bc^\top_y \by + \bc^\top_z \bz + \bc^\top_r \br + \sum_{w\in W} \theta_w \\
    \textrm{s.t.}         &\;    \theta_w \ge \; \phi^{(j)}_w + \bpi^{(j) \top}_w(\by - \by^{(j)})\notag \\
        &\qquad \qquad \;\;\,+ \blambda_w^{(j) \top}(\bz - \bz^{(j)}) \notag\\ 
        &\qquad \qquad  \;\;\,+ \bmu^{(j)\top}(\br - \br^{(j)}), \notag \\
        &\qquad \qquad \quad\;\forall j \in \{0, \dots, k-1\}, \; w \in W \
    \label{eq:CEM_cuts} \\
    &\; 0 \le \br \le \overline{\br}, \; \; \;  \by \in \mathbb{Z}_+^{n_y}, \;\;\; \bz \in \{ 0, 1 \}^{n_z} \label{eq:CEM_master_int_bin}
\end{align}
\end{subequations}
\noindent Here, $\theta_w$ are ``cost-to-go'' variables that provide a lower-bound estimate of the optimal cost of their respective operational sub-period.
The upper bound at iteration $k$ is given by:
\begin{align*}
    U^{(k)} = \min_{j = 0,...,k} \bc^\top_y \by^{(j)} + \bc^\top_z \bz^{(j)} + \bc^\top_r \br^{(j)} + \sum_{w \in W} \beta_w \bc^\top_w \bx_w^{(j)}
\end{align*}
Cutting planes \eqref{eq:CEM_cuts} are added after each iteration to improve the approximation of each sub-period's optimal value, $\theta_w$, as a function of investment decisions.

The BD algorithm is given in Algorithm \ref{alg:GBD}, following Pecci and Jenkins \cite{pecci2025regularized}. The algorithm first solves \eqref{eq:CEM_master} to obtain initial values for $\by^0$, $\bz^0$, and $\br^0$. It then solves the subproblem(s), \eqref{eq:CEM_subproblem}, to generate $U^{(k)}$ cuts, which are added to the master problem, and resolves the master problem to get $L^{(k)}$. This loop is repeated iteratively. If the values of $U^{(k)}$ and $L^{(k)}$ are within tolerance $\epsilon_{tol}$, the algorithm terminates and returns the best values found for $\by^*$, $\bz^*$, and $\br^*$.
\begin{algorithm}
    \caption{BD}\label{alg:GBD}
    \begin{algorithmic}[1]
    \State \textbf{Input:} Set $K_{max}$, $\epsilon_{tol}$, $k = 1$
    \State \textbf{Output:} $\by^*$, $\bz^*$, $\br^*$
    \State \;  \textbf{Initialization:} Solve \eqref{eq:CEM_master} with no cuts to obtain $\by^0$, $\bz^0$, and $\br^0$ 
    \While{$k \le K_{max}$}
    \For{ $w \in W$}
        \State Solve operational subproblem \eqref{eq:CEM_subproblem}
    \EndFor
    \State Compute upper bound $U^{(k)}$ and update $\by^*$, $\bz^*$, $\br^*$
    \State Update cuts in \eqref{eq:CEM_master} and solve to get $L^{(k)}$
    \If{$(U^{(k)} - L^{(k)}) / |L^{(k)}| > \epsilon_{tol}$}
        \State Stop and return $\by^*$, $\bz^*$, and $\br^*$
    \Else \; Set $\by^{(k + 1)}$, $\bz^{(k + 1)}$, $\br^{(k + 1)}$ to solution of \eqref{eq:CEM_master} \label{line:regularize_change}
    \EndIf
    \EndWhile
    
    \end{algorithmic}
\end{algorithm}
\vspace{0.1in}

\subsection{Generalized Benders Decomposition}

GBD \cite{geoffrion1972generalized} is an extension of BD for nonlinear problems with an identical iterative algorithmic structure. As implemented here, GBD differs from BD in two ways. 

\begin{enumerate}
    \item In GBD, DCOPF operational constraints (\ref{eq:CEM_subproblem_NLP_equality}) and (\ref{eq:CEM_subproblem_NLP_inequality}) take the form of (\ref{eq:DCOPF_bilinear}), resulting in a nonconvex, nonlinear program.
    \item  As a result, in (\ref{eq:CEM_subproblem}) and (\ref{eq:CEM_master}), dual variables $\bpi^{(k)}_w$, $\blambda^{(k)}_w$, and $\bmu^{(k)}_w$ are the Lagrange multipliers found from the Lagrangian of the bilinear problem rather than from the LP dual problem.
\end{enumerate}
\noindent The remainder of the algorithm is identical to BD.

Despite the change from dual variables to Lagrange multipliers in GBD, the form of cutting planes can be identical between BD and GBD (see, for instance, \cite{fischetti2017redesigning, mitrai2022multicut, pecci2025regularized}). Assuming full recourse between master and subproblems, each iteration of GBD results in a feasible solution since the master problem solutions are fixed in the subproblem. While BD and GBD use a different representation of the DCOPF constraints, both will result in the same primal solution (assuming there is a unique optimum). However, there is a small but important difference in the meaning behind the dual information obtained by BD and GBD. BD is applied to linear and mixed integer linear programs and is grounded in linear duality theory, where the values of $\bpi_w^{(k)}$, $\blambda^{(k)}_w$, and $\bmu^{(k)}_w$ would be dual values that can be found from the dual of \eqref{eq:CEM_subproblem}. In GBD, cuts are derived from nonlinear duality theory, and $\bpi_w^{(k)}$, $\blambda^{(k)}_w$, and $\bmu^{(k)}_w$ are the optimal Lagrange multipliers of the Lagrangian of \eqref{eq:CEM_subproblem}. This allows the use of \eqref{eq:DCOPF_bilinear} for the DCOPF constraints, which can result in different values in the cut constraint than would be returned using dual values in BD, potentially forming stronger cuts and recovering better solutions. 

Using \eqref{eq:DCOPF_bilinear} results in a nonconvex problem due to its bilinear constraints. This can result in the problem converging 
to a point that is not even a local optimum (see \cite{sahinidis1991convergence} and \cite{bagajewicz1991generalized}), so forming a valid lower bound on this problem is an additional point of interest that is discussed in Section \ref{sec:case_study}. Thus when GBD is implemented in Algorithm \ref{alg:GBD}, $U^{(k)}$ remains a valid upper bound since it is a feasible solution, but $L^{(k)}$ is not guaranteed to be a valid lower bound. Consequently, we will refer to $L^{(k)}$ as the ``algorithm bound'' in the case study below, and the termination criteria of $\epsilon_{tol}$ is not a true mixed integer gap. We will discuss in the next subsection how a valid mixed integer gap can be recovered.

\subsection{Strategies for Improving Performance}
In this section, we discuss three approaches to improving the performance of BD and GBD: hot-starting with a transport model, generating initial cuts from an LP relaxation of the master problem, and using a regularization scheme.

{\it Hot-Starting with Transport Model - } Romero and Monticelli \cite{romero2002hierarchical} propose hot-starting the algorithm by solving a transport-constrained (``pipe-flow'') version of the subproblems, then maintaining those cuts when DCOPF constraints are reintroduced. The transport model is equivalent to \eqref{eq:DCOPF_bilinear} without constraints \eqref{eq:DCOPF_exist}, \eqref{eq:DCOPF_bilinear_cand}, and \eqref{eq:DCOPF_bilinear_replace}. By construction, this is a relaxation of the subproblem constraints. Because the transport model removes all nonlinear equations, the resulting cuts provide a valid lower bound. Since the transport model is a relaxation, the lower bound remains valid for the original DCOPF-constrained GTEP problem. Thus, \eqref{eq:CEM} can be ``hot-started'' by first solving the transport model to optimality and maintaining those cuts in the master problem after DCOPF constraints are added. In this work, we will test first solving the transport-constrained problem with BD, maintaining those cuts, and then enforcing DCOPF constraints (e.g., \eqref{eq:DCOPF_exist}, \eqref{eq:DCOPF_bilinear_cand}, \eqref{eq:DCOPF_bilinear_replace}) and solving with BD or GBD.

A further benefit of solving the transport model is that it provides a valid lower bound on GBD, since GBD is not guaranteed to converge to a global or local optimum for nonconvex problems. Consequently, while $L^{(k)}$ provided by GBD can be used as a convergence criterion, it does not provide a valid optimality gap. However, if the transport model is solved, it constitutes a guaranteed lower bound and a valid optimality gap can be reported, and this could be used as an alternative termination criterion.

{\it LP Relaxations in Master Problem} To accelerate solution time, the computationally expensive integer decisions in the master problem (e.g. asset expansion or retirement decisions) can be relaxed to continuous decisions. This relaxed problem can be solved to optimality, and the resulting cuts maintained when integer constraints are reintroduced \cite{mcdaniel1977modified}. A modified version of this approach, implemented by Lumbreras et al. \cite{lumbreras2016solve}, relaxes all integer constraints in the master problem, then reimposes subsets of the original integer constraints until all are reinforced. In this work, we will test solving an LP relaxation of the master problem to optimality with GBD, maintaining the generated cuts, and then solving the integer-constrained master problem.

{\it Regularization -} In both BD and GBD, the master problem can choose extreme solutions that elicit poor cuts from the subproblems. Regularization (also referred to as stabilization) modifies the master problem to choose sub-optimal, less extreme solutions, resulting in stronger elicited cuts \cite{lemarechal1995new, pecci2025regularized, ruszczynski1986regularized}. We regularize our problem via the approach in Pecci and Jenkins \cite{pecci2025regularized}. Their approach solves a regularized master problem defined by \eqref{eq:CEM_master_reg}.
\begin{subequations}\label{eq:CEM_master_reg}
\begin{align}
    \min &\; \Phi(\by, \bz, \br) \\
    \textrm{s.t.} &\; \theta_w \ge \; \phi^{(j)}_w + \bpi^{(j) \top}_w(\by - \by^{(j)}) \notag\\ 
        & \qquad \qquad \;\;\,+ \blambda_w^{(j) \top}(\bz - \bz^{(j)})\notag\\
        & \qquad \qquad \;\;\,+ \bmu^{(j)\top}(\br - \br^{(j)}),\notag\\
        &\qquad \qquad \qquad\qquad\,\forall j \in \{0, \dots, k-1\}, \; w \in W\\
    &\; L_\alpha^{(k)}\ge \bc^\top_y \by + \bc^\top_z \bz + \bc^\top_r \br + \sum_{w\in W} \theta_w  \label{eq:CEM_master_reg_obj_bound}\\
    &\; 0 \le \br \le \overline{\br}, \; \; \;  \by \in \mathbb{Z}_+^{n_y}, \;\;\; \bz \in \{ 0, 1 \}^{n_z}
\end{align}
\end{subequations}
Here, $L_\alpha^{(k)} = L^{(k)} + \alpha (U^{(k)} - L^{(k)})$ where $\alpha$ is a scalar value between 0 and 1. The function $\Phi(\cdot)$ can take various forms; we use $\Phi(\by, \bz, \br) = 0$. \pref{eq:CEM_master_reg} is, then, a feasibility problem where \eqref{eq:CEM_master_reg_obj_bound} requires the solution be within a bound, set by a fraction of the optimality gap, of the optimal objective value of \eqref{eq:CEM_master}. The resulting solution  can be more interior when $\alpha$ is large, and closer to optimal when small. The regularized algorithm is almost identical to Algorithm \ref{alg:GBD}, but Line \ref{line:regularize_change} is adjusted to solve and return solutions from \eqref{eq:CEM_master_reg} instead of returning solutions from \eqref{eq:CEM_master}.

In Section \ref{sec:case_study}, we will test a combination of these three approaches. Algorithm \ref{alg:GBD_improved} shows the procedure with transport model hot-starting and LP relaxation. There are four steps: 1) solve the LP relaxation of the GTEP model with transport constraints, 2) solve the MILP form of the GTEP model with transport constraints, 3) solve the LP relaxation of the GTEP model with DCOPF constraints, and 4) solve the MILP form of the GTEP model with DCOPF constraints. The cuts obtained at each step are maintained in all succeeding steps to ``hot-start'' the next step. 
\begin{algorithm}
    \caption{BD/GBD with Speedup Strategies}\label{alg:GBD_improved}
    \begin{algorithmic}
    \State \textbf{Input:} Set $K_{max}$, $\epsilon_{tol}$, $k = 1$
    \State \textbf{Output:} $\by^*$, $\bz^*$, $\br^*$
    \State \;  \textbf{Initialization:} Solve \eqref{eq:CEM_master} with no cuts to obtain $\by^0$, $\bz^0$, and $\br^0$ 
    \State \textbf{Step 1: LP Relaxation of Transport Model} - Relax integrality constraints \eqref{eq:CEM_master_int_bin}. Apply Algorithm \ref{alg:GBD} using the transport model for \eqref{eq:CEM_subproblem_NLP_equality} and \eqref{eq:CEM_subproblem_NLP_inequality}.
    \State \textbf{Step 2: Integer Form of Transport Model} - Maintain cuts in \eqref{eq:CEM_master} from Step 1. Add integrality constraint of \eqref{eq:CEM_master_int_bin}. Apply Algorithm \ref{alg:GBD} using the transport model for \eqref{eq:CEM_subproblem_NLP_equality} and \eqref{eq:CEM_subproblem_NLP_inequality}.
    \State \textbf{Step 3: LP Relaxation of DCOPF Model} - Maintain cuts in \eqref{eq:CEM_master} from Steps 1 and 2. Relax integrality constraints \eqref{eq:CEM_master_int_bin}. Apply Algorithm \ref{alg:GBD} using DCOPF constraints \eqref{eq:DCOPF_linear} or \eqref{eq:DCOPF_bilinear} for \eqref{eq:CEM_subproblem_NLP_equality} and \eqref{eq:CEM_subproblem_NLP_inequality}.
    \State \textbf{Step 4: Integer Form of DCOPF Model} - Maintain cuts in \eqref{eq:CEM_master} from Steps 1, 2, and 3. Add integrality constraint of \eqref{eq:CEM_master_int_bin}. Apply Algorithm \ref{alg:GBD} using DCOPF constraints \eqref{eq:DCOPF_linear} or \eqref{eq:DCOPF_bilinear} for \eqref{eq:CEM_subproblem_NLP_equality} and \eqref{eq:CEM_subproblem_NLP_inequality}.    
    \end{algorithmic}
\end{algorithm}
\vspace{0.1in}

 Importantly, any of the 4 steps could be omitted yet maintain a valid algorithm, and regularization could be applied to any step. In this work, we test combinations of Algorithm \ref{alg:GBD_improved} as outlined in Table \ref{tab:strategies} along with the abbreviation used to notate the results.

\begin{table}\centering
\scriptsize
\caption{Speedup Strategies - {\it HS} Refers to hot-starting with the transport model, {\it LP} refers to using an LP relaxation, {\it Reg} refers to using regularization, and {\it SemiReg} refers to using regularization for only some steps}\vspace{4pt}
\label{tab:strategies}
\begin{tabular}{lll}
\hline
   Name & Abbreviation & Details \\
\noalign{\smallskip}\hline\noalign{\smallskip}
    Baseline & - &Algorithm \ref{alg:GBD}\\
    Strategy 1 & HS & Steps 2 and 4 of Algorithm \ref{alg:GBD_improved}  \\
    Strategy 1 & LP & Steps 3 and 4 of Algorithm \ref{alg:GBD_improved} \\
    Strategy 1 & Reg & Step 4 of Algorithm \ref{alg:GBD_improved} with regularization \\
    Strategy 2 & HS + LP & Steps 1-4 of Algorithm \ref{alg:GBD_improved} \\
    Strategy 3 & HS + LP + Reg & Steps 1-4 of Algorithm \ref{alg:GBD_improved} with regularization  \\
    Strategy 4 & HS + LP + SemiReg & Steps 1-3 of Algorithm \ref{alg:GBD_improved} with regularization, \\&&then step 4 without regularization\\
\noalign{\smallskip}\hline
\end{tabular}\vspace{-15pt}
\end{table}

\section{Case Study} \label{sec:case_study}

\subsection{System Overview}
\begin{figure}[!t]
    \centering
    \includegraphics[width=1\linewidth]{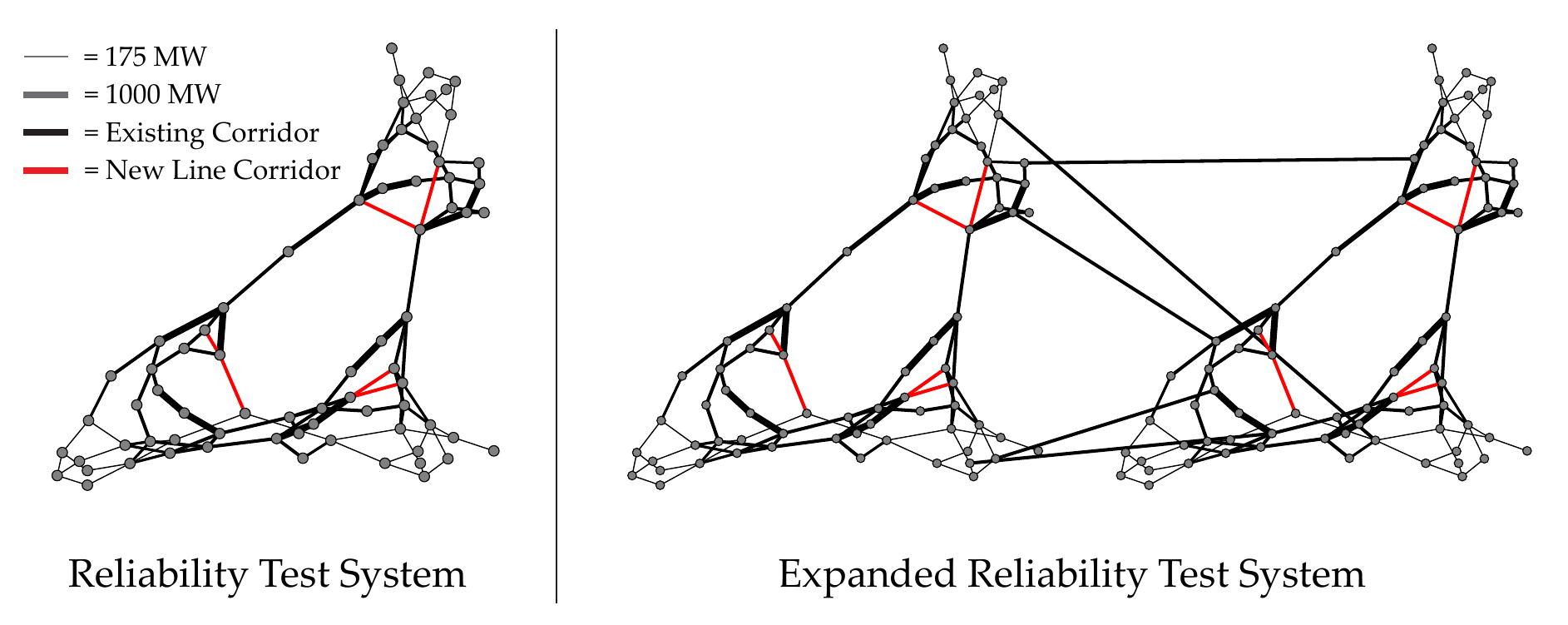}
    \caption{A visualization of the RTS and expanded RTS. Black lines indicate existing corridors while red lines indicate new corridors, with the size of the line corresponding to the capacity of the line(s) on the corridor.}
    \label{fig:RTS_systems}
\end{figure}

In this case study, we solve three systems. The systems we solve are adaptations of the Reliability Test System (RTS) \cite{barrows2019ieee}. The RTS has 73 buses split into three zones, with 108 transmission corridors and 120 existing lines. 
We add an additional 6 corridors where a new line can be constructed. To test algorithm scalability, we construct an expanded RTS by duplicating the original RTS, perturbing load and generator variability data in this copy, and arbitrarily adding five transmission lines between the systems. A single zone of the system will also be solved as a smaller test case. The system variants are visualized in Figure \ref{fig:RTS_systems}. The original RTS contains 5.42 GW of renewable power and 9.45 GW of thermal power with hourly data for one year and peak load of 8.91 GW. Because the original RTS has surplus transmission and generation, we scale the load by 2 to require new generation and transmission construction.

The original RTS does not contain candidate line information, so synthetic candidate information had to be generated. For each existing line, we allow a new ``candidate line’’ based on the existing line's data and location. Existing lines over 200 MW were treated as fixed (i.e., they belong to $\mathcal{L}_{fixed}$). A candidate line was added parallel to each of these lines with a capacity 1.5 times the existing line, and a reactance 1.5 times smaller than the existing line. In addition, lines in $\mathcal{L}_{fixed}$ could be reconductored. Reconductoring is treated in a piecewise manner: it can increase capacity up to 10\% for a lower cost than constructing a new line and can increase an additional 15\% capacity for a higher cost. This represents the real-world ability to increase capacity in existing lines by minor interventions like tightening cables or major interventions like replacing cables with new materials. Existing lines less than 200 MW were treated as lines that could be replaced (i.e., they belong to $\mathcal{L}_{rep}$). Each corresponding candidate replacement line was 2.5 times the capacity with a reactance 2.5 times smaller than the original line. This process resulted in 120 new candidate lines, with $|\mathcal{L}_{fixed}| = 83$ and $|\mathcal{L}_{rep}| = 37$. Six new corridors were also added, 2 in each zone. These new corridors were determined by placing lines between the bus with the largest demand and the two buses with the largest generating capacity not originally directly linked to the largest demand bus. The value of $M_\ell$ in \pref{eq:DCOPF_linear} was set to the smallest value we could determine for existing corridors. For existing corridors without reconductoring, this value was $\overline{f}_\ell$, with  $ 1.25  \overline{f}_\ell$ on the other existing corridors with reconductoring. For new corridors, we set $M_\ell$ to $10  \overline{f}_\ell$ 
to ensure phase angles are unrestricted. 

Length and price data for new lines also had to be generated, as they are not in the original RTS data. For new transmission, we used the cost of \$3,667 per MW-mile reported by \cite{Ho2021} for spur lines. We converted these values from 2015 dollars to 2023 dollars by multiplying by 1.233 (based on inflation index values obtained in the Regional Energy Deployment System source code). We used line lengths reported by \cite{grigg1999ieee}. As some line lengths were listed as zero, we set a minimum annuitized line cost of \$20,358 (cost of a 5 mile line length). As line lengths were almost identical for each zone, we added random noise of up to 10\% of the total cost to help avoid degeneracy in the final solutions. Annuitized costs for each line were generated using a weighted average cost of capital (WACC) of 4.4\% and a 60 year asset life.

Ten different new generation technologies were available within each zone with each technology having three locations it could be placed within each zone (see documentation accompanying our linked code for details on generator placement). We require each to be built in 200 MW increments to prevent unrealistically small capacity decisions. 



We used a time domain reduction (TDR) strategy to reduce the RTS data from 52 weeks to 16 representative weeks (i.e., $|W| = 16$). Our TDR includes the extreme time periods (minimum solar, minimum wind, maximum demand) then uses k-means clustering to group the rest of the year into representative weeks (see GenX.jl documentation \cite{GenXsoftware}). 

Code and data for reproducing the following results can be found at \cite{Zenodo_data}. All problems were solved with Gurobi on an AMD EPYC 9654 2.4 GHz processor using 32 cores and 160 and 192 GB of RAM for the RTS and expanded RTS, respectively. The MIP master problem was solved to a 0.1\% gap or to a time limit of 6 min (if the problem timed out, the best incumbent was used). The quadratically constrained subproblems were likewise solved, and the dual retrieved with Gurobi. Solving the dual KKT system inside Gurobi can be numerically unstable, such that Gurobi does not always return a dual value. When this occurred, the affected cut was skipped.
The single zone, RTS, and expanded RTS problems had 1,645,179 (40 integer, 30 binary), 5,037,695 (90 integer, 126 binary), and 10,102,269 (180 integer, 252 binary) variables, with 2,483,959, 7,613,328, and 15,264,288 constraints, respectively. 

\subsection{Results and Discussion}\label{Sec:Results}

\subsubsection{Baseline and Big-M Tuning}
\begin{figure}[!t]
    \centering
    \includegraphics[width=1\linewidth]{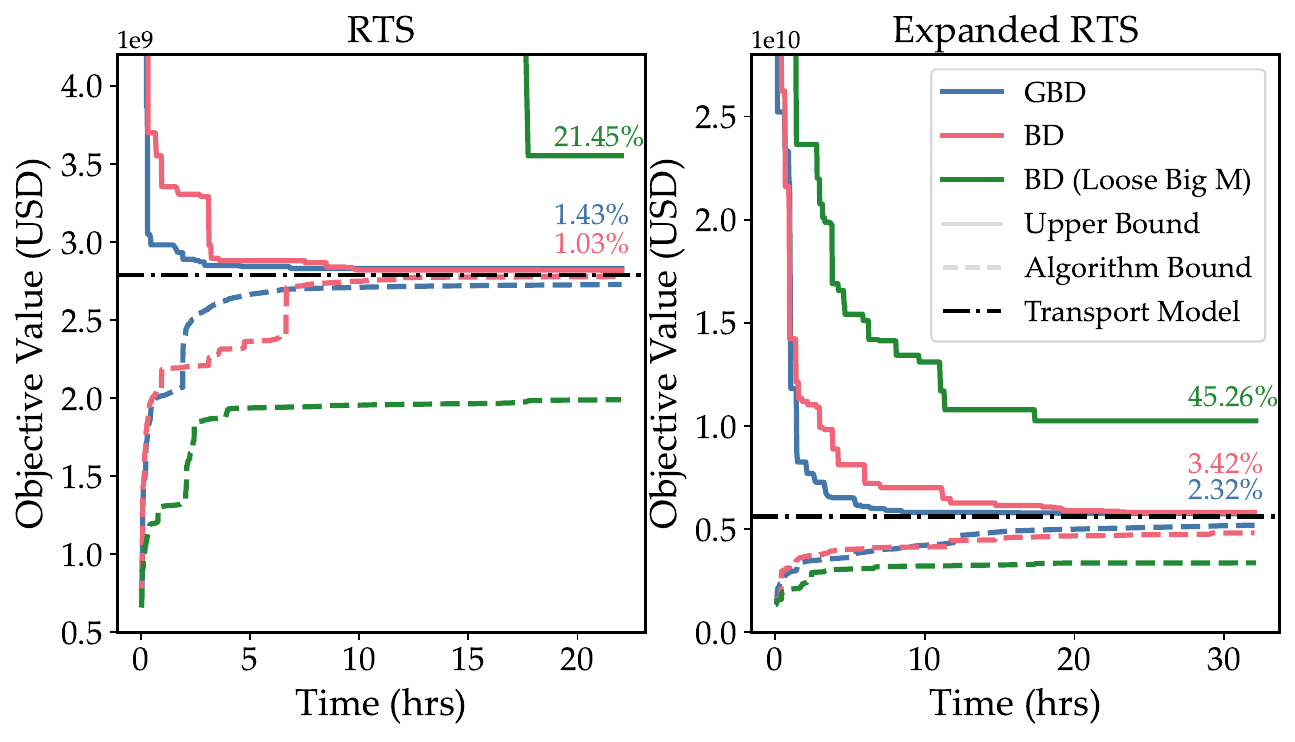}\vspace{4pt}
    \caption{Results of GBD and BD applied to the RTS and expanded RTS with no speedup strategies. The solution of a transport constrained model is shown in black, with the gap between the best solution of the algorithm and the transport constrained model annotated on the right-hand side. ``Loose Big M'' uses a value of $M_\ell = 10 \overline{f}$. }
    \label{fig:baseline}
\end{figure}

As a baseline, Figure \ref{fig:baseline} shows the performance of BD and GBD with no speed-up strategies on both the RTS and expanded RTS. To highlight the importance of tuning $M_\ell$ in \eqref{eq:DCOPF_linear}, we also show the performance of BD with $M_\ell = 10 \overline{f}_\ell, \ell \in \ml$, here referred to as ``loose Big M''. For the RTS, BD ultimately recovered the best overall solution by a margin of over \$10 million. For the expanded RTS, GBD outperformed BD by more than \$65 million. In both cases, the ``loose Big M'' solution performed poorly, never getting within 20\% of the transport constrained solution. However, all of these approaches outperformed a monolithic approach. Solving either the RTS or the expanded RTS as a monolithic problem with Gurobi with 32 cores was intractable with no incumbent returned in 24 and 32 hours, respectively. 

\subsubsection{Speedup Strategy Comparison}
\begin{figure}[!t]
    \centering
    \includegraphics[width=0.95\linewidth]{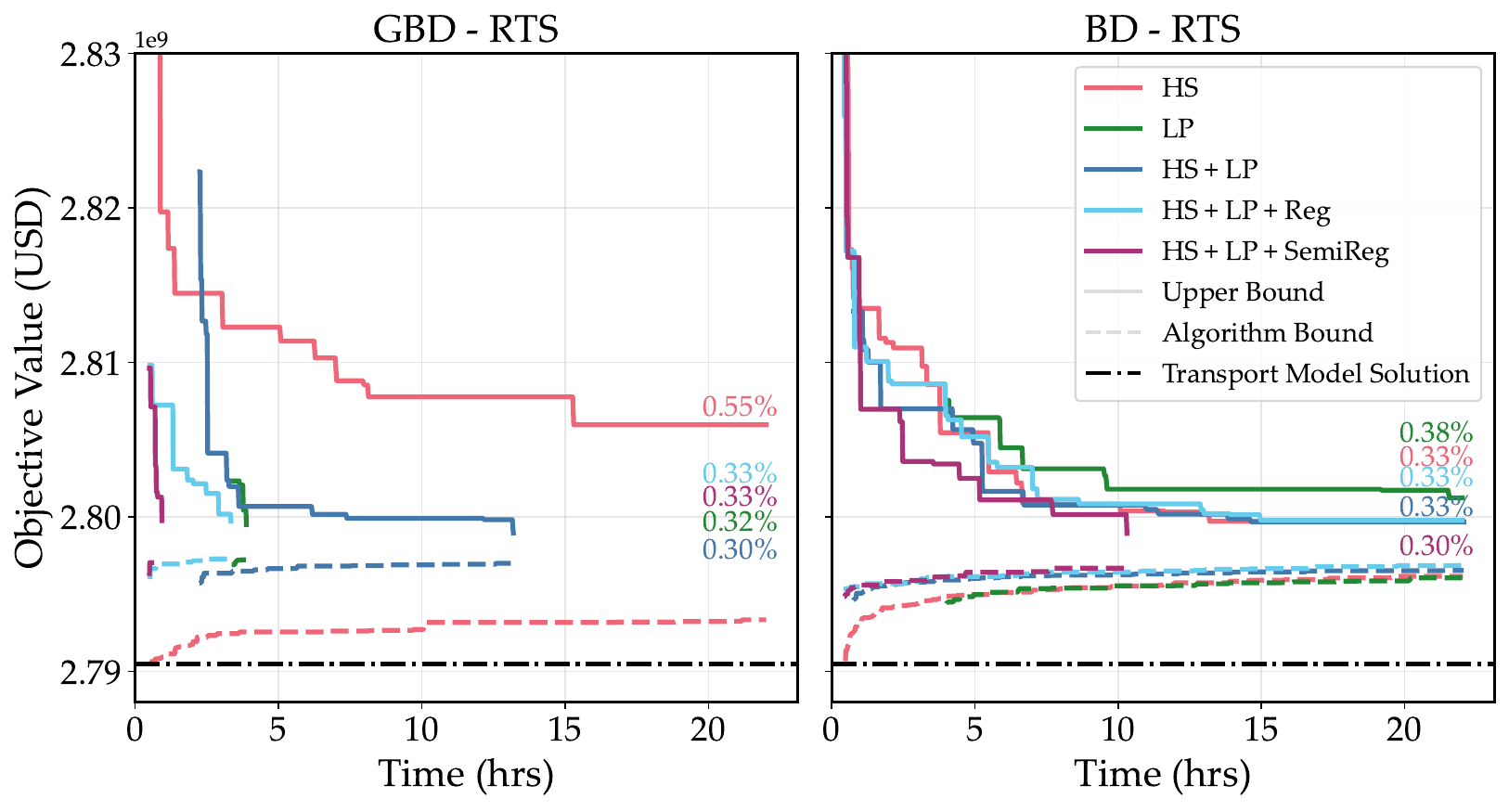}
    \caption{Comparison of the six speedup strategies (Table \ref{tab:strategies}) for GBD and BD on the RTS. Gaps reported on the right-hand side are compared to the solution of the transport constrained model. ``Reg'' is not shown because it performed worse than alternatives ($>$1\% gap). Lines begin after any hot-starting and LP relaxation procedures finish.}
    \label{fig:GBDvBD_RTS}
\end{figure}

\begin{figure*}[!t]
    \centering
    \includegraphics[width=.8\linewidth]{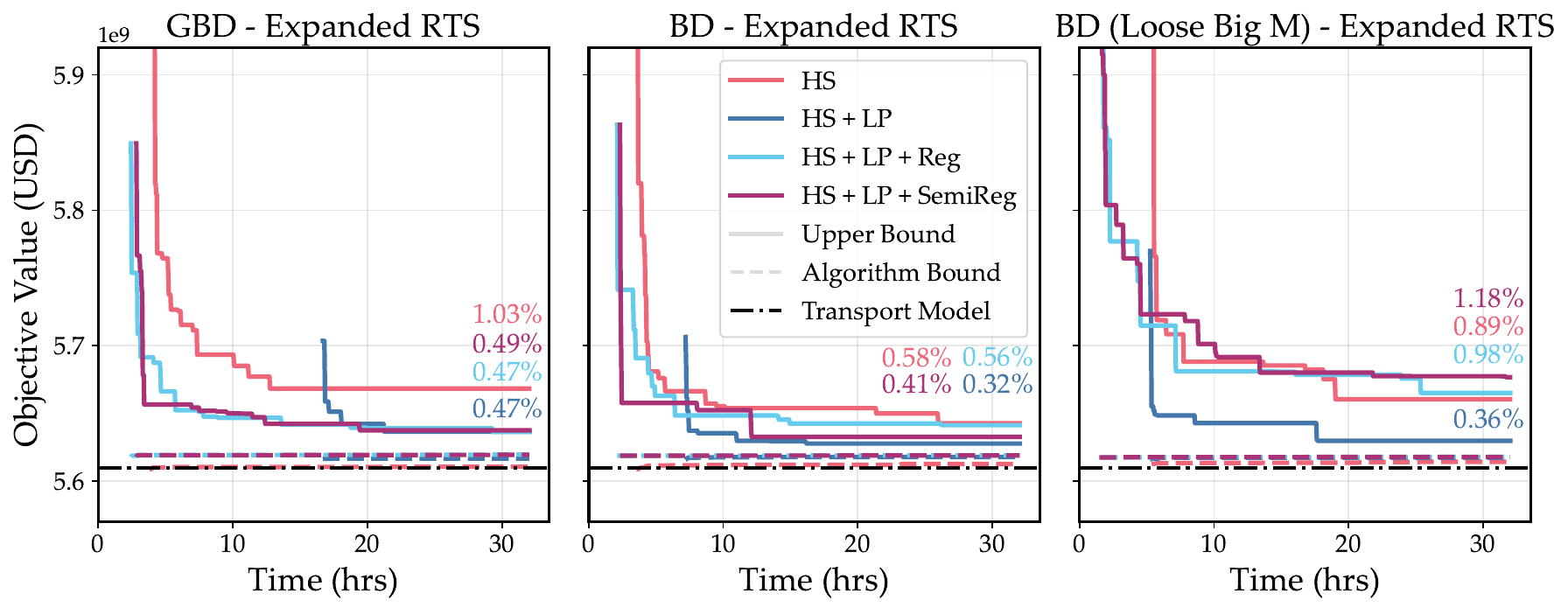}
    \caption{Comparison of the six speedup strategies (Table \ref{tab:strategies}) for GBD and BD on the expanded RTS. Gaps reported on the right-hand side are compared to the solution of the transport constrained model. ``LP'' and ``Reg'' are not shown because the LP relaxation did not finish in the given time, and the regularization alone never got below a 5\% gap. Lines begin after any hot-starting and LP relaxation procedures finish.} 
    \label{fig:GBDvBD_expanded_RTS}
\end{figure*}

The gaps shown in Figure \ref{fig:baseline} can be further reduced by utilizing the speedup strategies from Table \ref{tab:strategies}. The results of the speedup strategies for the RTS and expanded RTS are shown in Figures \ref{fig:GBDvBD_RTS} and \ref{fig:GBDvBD_expanded_RTS}, respectively with reported gaps based on the solution of the transport model. The ``Reg'' strategy is not shown due to its poor performance (never below a 1\% and 5\% gap for the RTS and expanded RTS, respectively). For the RTS, the speedup strategies ``HS'' and ``LP'' both significantly improved performance, while ``Reg'' had little impact and failed to decrease the overall solution below a 1\% gap. For GBD, ``HS'' performed the next worst, and all other strategies reached a comparable final solution and terminated when at $\epsilon = 0.001$. For these, the regularization scheme helped the overall problem converge faster, with the ``HS + LP + SemiReg'' ultimately converging the fastest, using under 2 hours. For BD, the ``LP'' strategy performed the worst after ``Reg'', with ``HS'', ``HS + LP'', ``HS + LP + Reg'', and ``HS + LP + SemiReg'' reaching comparable solutions, with only ``HS + LP + SemiReg'' reaching the termination criteria of $\epsilon = 0.001$.  

Results for the expanded RTS are given in Figure \ref{fig:GBDvBD_expanded_RTS}. For GBD, the ``HS + LP'', ``HS + LP + Reg'', and ``HS + LP + SemiReg'' strategies all performed comparably, taking over 20 hours to get below a 0.5\% gap. For BD, the ``HS + LP'' and ``HS + LP + SemiReg'' strategies performed better than for GBD, obtaining slightly better solutions in less overall time. For the loose Big M case, only the ``HS + LP'' strategy recovered a solution under a 0.5\% gap. Regularization does significantly speedup the initial steps for both GBD and BD. For GBD and BD, the ``HS + LP'' strategy required more than 15 and 5 hours, respectively, to finish the initial hot-starting and LP relaxation. However, for the case with the looser big M value, the problems using regularization failed to recover as high quality solutions as the other methods, suggesting regularization is more sensitive to the value of big M. Unlike the RTS case, no strategy with GBD or BD converged to the 0.001 convergence tolerance.  

\begin{figure*}[!t]
    \centering
    \includegraphics[width=0.8\linewidth]{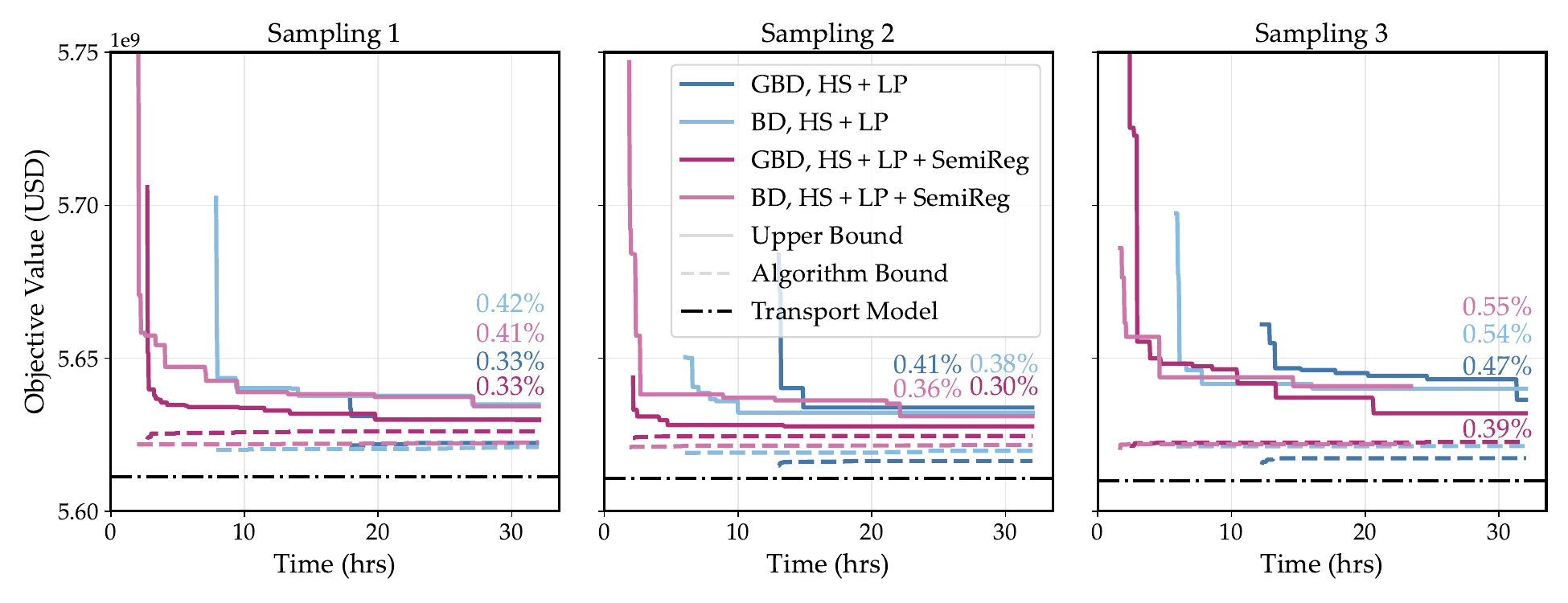}
    \caption{Results of the expanded RTS system with three different random samplings of available lines for new line construction. }
    \label{fig:RTS_samplings}
\end{figure*}

While the results shown in Figures \ref{fig:GBDvBD_RTS} and \ref{fig:GBDvBD_expanded_RTS} allow for new line constructions on any existing corridor, real systems will not have such flexibility. To more closely represent a real system, we solve three instances of the expanded RTS with fewer available candidate lines and reconductorable lines by taking three random samplings of 80 lines from $\mathcal{L}_{cand}$ and 20 lines from $\mathcal{L}_{rec}$ and solving the model. For simplicity, we use the ``HS + LP'' and ``HS + LP + SemiReg'' strategies because they performed well for both GBD and BD in the tests presented above. The results of these three systems with random line samplings are shown in Figure \ref{fig:RTS_samplings}. In all three cases, GBD and BD are both able to achieve high-quality solutions, but GBD with strategy ``HS + LP + SemiReg'' performed the best in all three cases. In samplings 1 and 3, GBD ultimately recovered the best overall solutions compared with BD.

\begin{figure}[!t]
    \centering
    \includegraphics[width=1\linewidth]{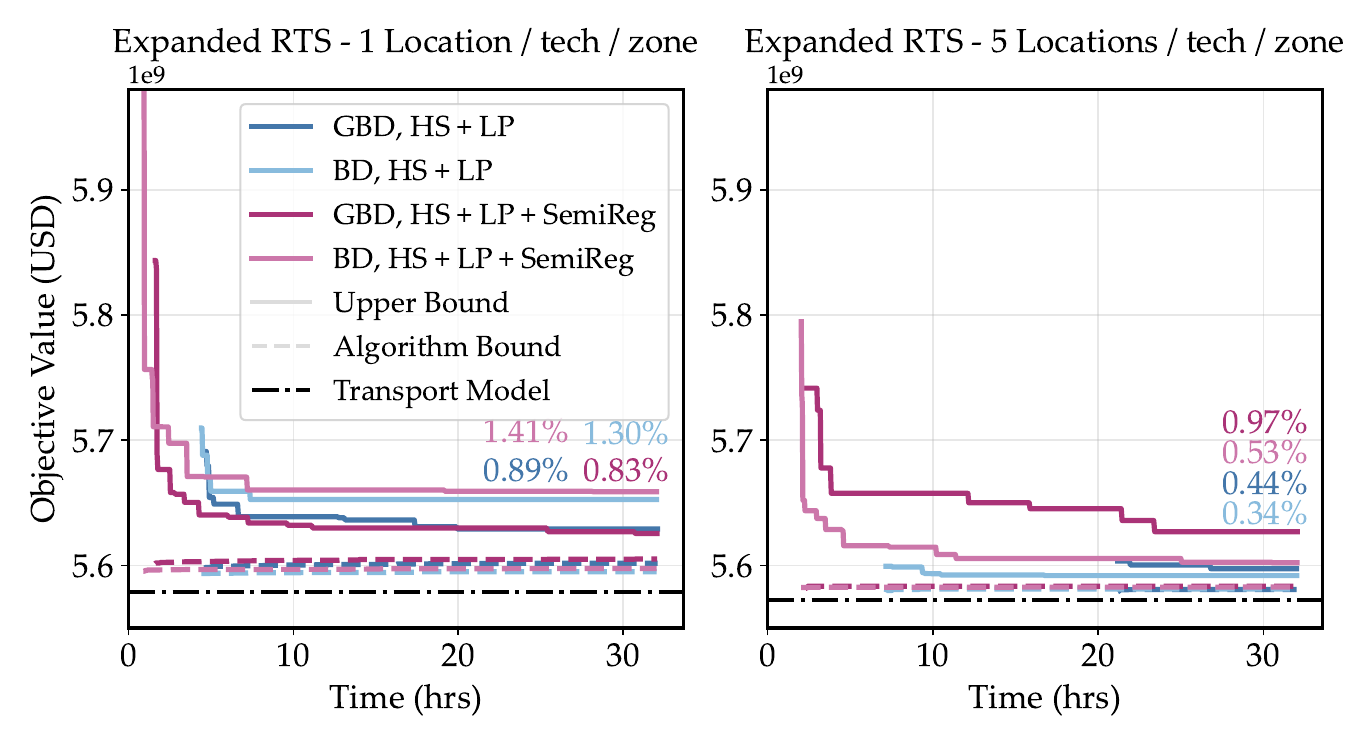}
    \caption{Results for GBD and BD with two different speedup strategies for the expanded RTS with 1 and 5 locations per technology per zone available. Reported gaps are from the optimal solution of the transport model.}
    \label{fig:RTS_gen_locations}
\end{figure}

As an additional test, we changed the number of generator locations available in each zone and tested with 1 and 5 locations per technology per zone. In the above cases, each technology could be built on three buses in each zone (see case study documentation). For the 1 location per technology per zone case, the one location was a subset of the three locations used for that technology in the previous analysis. For 5 locations per technology per zone, we added two additional locations to the three locations used. The 1 and 5 locations per technology per zone solutions for the expanded RTS are shown in Figure \ref{fig:RTS_gen_locations}. GBD performs noticeably better than BD for the 1 location per technology per zone case, though it never gets below a 0.5\% gap. In the 5 location per technology per zone case, BD also performed better than GBD, with the ``HS + LP'' strategy achieving the best solutions. We hypothesize that GBD performs better in the 1 location per technology per zone case because it does not use the ``big M'' formulation which may help it produce stronger cuts for transmission expansion problems and thus better identify the value of new transmission lines. This hypothesis explains the performance trend because increasing the number of possible locations per tech per zone could reduce the need for new transmission or introduce degenerate solutions in new line constructions as more modular generation can be built. This could also explain why we obtained smaller gaps in the 5 location case: the transport model solution could be closer to the DCOPF constrained model as the model can substitute more generation for new transmission. Confirming this hypothesis would require future analysis. 


\subsubsection{Scalability Comparison}
\begin{figure}[!t]
    \centering
    \includegraphics[width=0.6\linewidth]{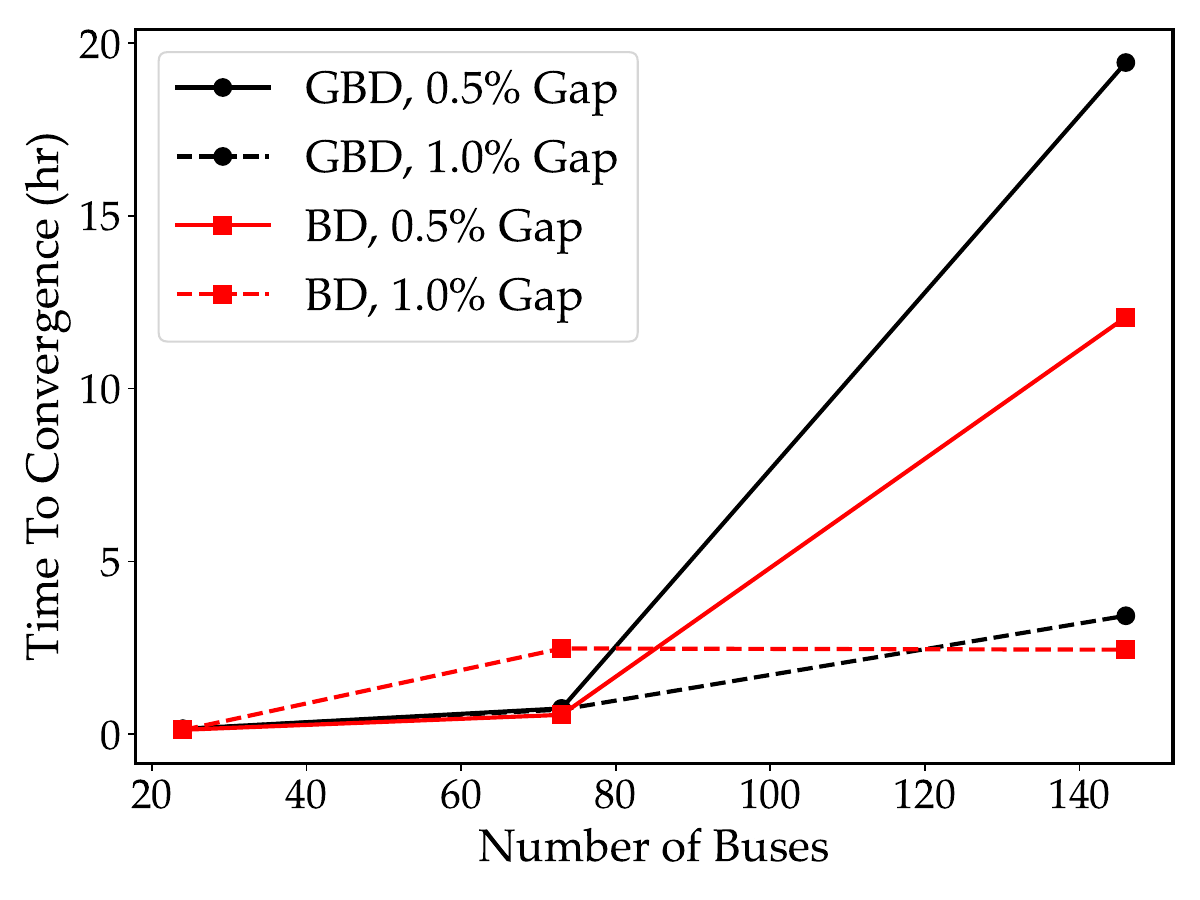}
    \caption{Time to convergence to a 1\% and 0.5\% gap (compared to the solution of the transport model) as a function of the number of buses in the system for GBD and BD with the ``HS + LP + SemiReg'' strategy.}
    \label{fig:GBD_scaling}
\end{figure}

Finally, to test the scalability of GBD, we also solved a single zone of the RTS containing 24 buses. With BD and ``HS + LP + SemiReg'', the single zone problem is solved in 454 seconds and 63 iterations to a tolerance of 0.1\% (compared with the algorithm bound), while GBD with the ``HS + LP + SemiReg'' required 551 seconds and 31 iterations. For the smaller system, BD is significantly faster per iteration, though both BD and GBD reach solutions efficiently. This additional solve time is due in part to solving the quadratically constrained subproblems of GBD, which is slower than solving the linear relaxation. The solve times of the three system sizes are compared in Figure \ref{fig:GBD_scaling} for GBD and BD with the ``HS + LP + SemiReg'' strategy for both a 0.5\% and 1.0\% gap (compared with the transport constrained model; the reported time for the 1 zone, 24 bus case is to convergence since the final solutions were greater than a 0.5\% gap with the transport model). The performance of both BD and GBD to obtain a 0.5\% gap scales nonlinearly with the number of buses. A key challenge with GTEP problems is solving large-scale systems. While decomposition provides a powerful tool to make these problems tractable, a system on the scale of several hundred or thousands of buses will likely be intractable with GBD or BD. Importantly, much of the time in both BD and GBD is spent shrinking the gap after a reasonably high-quality solution is reached. This is seen in Figure \ref{fig:GBD_scaling}, where a 1.0\% gap is obtained faster and scales better than going to a smaller gap. Interestingly, BD actually achieved a sub 1\% gap with the expanded RTS faster than the smaller RTS case. This highlights some of the nature of BD or GBD in that there is some randomness in the path that the algorithm takes to reach convergence.

\subsubsection{Study Limitations and Assumptions}
This study does include some limitations and assumptions. First, results may be system dependent, and extrapolating to other systems with different characteristics (e.g., higher voltage lines, different generation mixes, different levels of congestion) could impact algorithm performance to some degree. The RTS is a synthetic system, though widely used, and thus naturally differs from many systems in use today. Further, the generation mix from one system to another can differ between regions or countries, and this generation mix could likewise change the amount of new generation or new transmission required in a system. Likewise, problem data (e.g., fuel prices, investment costs) impact system operation but change regularly in the real world due to political policies, world events, or technology development. Lastly, the BD algorithm can on occasion simply get ``lucky'' in its path to a high quality solution because there are no guarantees that each iteration of Benders improves on the feasible solution of the previous iteration. These limitations are a part of any empirical study like this work. Addressing these challenges was part of the motivation for testing with different line samples and different number of generator locations for new technologies, ensuring that we tested different system configurations and sizes. We worked to ensure that the data was realistic based on reliable databases. We likewise made certain that our compute resources were the same across each trial.

\section{Conclusion}\label{sec:conclusion}
Overall, we find the following conclusions: i) Decomposition strategies are essential to solve GTEP of any realistic scale. In our case, the RTS and expanded RTS were intractable without decomposition in the same time period. ii) Speedup strategies are essential for getting GTEP problems of even medium size to converge with BD or GBD.  iii) GBD is a useful strategy that in some cases (not all) can reach better solutions than the same methods with BD. iv) GBD alleviates the challenge of tuning $M_\ell$ in more common `Big M' relaxation formulations of the GTEP problem. In this work, the value of $M_\ell$ could be set at or close to $\overline{f}_\ell$ as the capacity and susceptance of the new candidate lines was directly proportional to existing lines. In other cases, this might not be possible, and other computational methods such as computing weighted shortest paths \cite{binato2002new} can be used for determining $M_\ell$. GBD does provide the benefit that this determination would not be required. If values of $M_\ell$ are not tuned, regularization strategies may be less effective. v) BD and GBD appear to require a combination of hot-starting with the transport model and a linear relaxation for best performance for GTEP problems. Overall, the ``HS + LP + SemiReg'' strategy seemed to perform best for the systems tested.  To our knowledge, no other studies solving GTEP problems have tested GBD compared to BD, and no other studies have analyzed these three combined speedup strategies for the performance of either algorithm. Further, no other GTEP studies to our knowledge consider the strategy of turning off regularization after a certain point to further improve solutions.

There are also several areas of future work to which this work could be applied. While solving the subproblem \eqref{eq:CEM_subproblem} using the quadratic constraints of \eqref{eq:DCOPF_bilinear} compared with the linear big M formulation of \eqref{eq:DCOPF_linear} returns the same primal solutions (assuming no degeneracy), solving with the quadratic constraints often resulted in longer subproblem solve times. If the Lagrange multipliers could be recovered from an LP representation of \eqref{eq:DCOPF_bilinear}, GBD could be sped up. In addition, since GBD does not guarantee a lower bound, the gaps reported in this work were based on the transport constrained solution. However, GBD often reported a tighter algorithm bound than BD. It is possible that a tighter global lower bound can be obtained for GBD problems. 

The systems addressed in this work were relatively small compared with real systems, which can contain hundreds to thousands of buses. A key challenge in this work was designing a test system that considers both transmission and generation expansion with realistic data. To our knowledge, larger systems are not readily available in easy-to-use data formats and would require additional work to create the required dataset for this analysis. It is likely that system size will eventually become computationally intractable as a mixed integer problem, and additional methods may need to be developed to recover high quality solutions. One option is to use further relaxations, such as only solving the LP solution of the overall model, which may be informative for many use cases. Such results are implicitly included in this work (e.g., for the ``HS + LP + SemiReg'' strategy) which suggest that LP relaxations could be tractable for larger-scale systems. In addition, there are further combinations of approaches that could be used to improve speedup, such as running BD with the ``HS + LP + SemiReg'' strategy, and then switching to GBD after some time. This hybrid approach could provide a stronger lower bound than the transport model and could result in better overall solutions.

\section*{Acknowledgements}
The authors thank Dr. Jose Daniel Lara and the Sienna team at the National Laboratory of the Rockies for providing the RTS data and technical support. They also thank Drs. Greg Schivley and Luca Bonaldo for their technical support.


\bibliographystyle{IEEEtran}
\bibliography{GBD_abb}

\end{document}